\documentclass{amsart}
\usepackage{graphicx}
\usepackage{amssymb}
\usepackage{graphicx}
\usepackage[all]{xy}

\theoremstyle{plain}

\newtheorem{theorem}{Theorem}[section]
\newtheorem{corollary}[theorem]{Corollary}

\newtheorem{lemma}[theorem]{Lemma}
\newtheorem{proposition}[theorem]{Proposition}

\theoremstyle{definition}
\newtheorem{definition}[theorem]{Definition}

\newtheorem{assumption}[theorem]{Assumption}

\theoremstyle{remark}

\newtheorem{remark}[theorem]{Remark}
\newtheorem*{remark*}{Remark}

\numberwithin{equation}{section}
\numberwithin{figure}{section}
\let\rom\textup
\def\a{\alpha}
\def\si{\sigma}
\def\ov{\overline}
\def\ciM{M^\circ}
\DeclareMathOperator\im{Im} \DeclareMathOperator\ind{ind}
 \DeclareMathOperator\Ell{Ell}

\newcommand{\Vect}{{\rm Vect}}
\begin{document}

\language0

\title[Elliptic Theory on Manifolds with Corners: II]{Elliptic Theory
on Manifolds with Corners:\\ II. Homotopy classification and $K$-Homology}

\author[Nazaikinskii]{Vladimir Nazaikinskii}

\address{%
Institute for Problems in Mechanics, Russian Academy of Sciences,
\newline\indent pr.~Vernadskogo 101-1, 119526 Moscow, Russia}

\email{nazaikinskii@yandex.ru}

\thanks{Supported in part by RFBR grants 05-01-00982
and 06-01-00098, by President of the Russian Federation grant
MK-1713.2005.1, and by the DFG project 436 RUS
113/849/0-1\protect\raisebox{1pt}{\circledR}\ ``K-theory and Noncommutative
Geometry of Stratified Manifolds"}
%----------Author 2
\author[Savin]{Anton Savin}

\address{Independent University of Moscow,
Bol'shoi Vlas'evskii per.~11,\newline\indent 119002 Moscow, Russia}

\email{antonsavin@mail.ru} \email{sternin@mail.ru}

%----------Author 3
\author[Sternin]{Boris Sternin}



% ----------------------------------------------------------------
\begin{abstract}
We establish the stable homotopy classification of elliptic pseudodifferential
operators on manifolds with corners and show that the set of elliptic operators
modulo stable homotopy is isomorphic to the $K$-homology group of some
stratified manifold. By way of application, generalizations of some recent
results due to Monthubert and Nistor are given.
\end{abstract}
\maketitle
\tableofcontents
% ----------------------------------------------------------------
\section*{Introduction}

Recently there has been considerable progress in understanding the notion of
ellipticity on noncompact manifolds and manifolds with singularities. For a
wide class of manifolds, ellipticity conditions for operators were established
and the corresponding finiteness theorems\footnote{Stating that an elliptic
operator is Fredholm in certain function spaces.} were proved; the
corresponding operator $C^*$-algebras were constructed. Hence the study of
\textit{topological} aspects of the theory of elliptic operators becomes
topical. Here one mainly speaks of the classification problem and the index
problem. Note that Gelfand's homotopy classification problem for elliptic
operators can naturally be restated in modern language as the problem of
computing the $K$-groups of symbol algebras, which prove noncommutative in most
cases. Thus Gelfand's problem naturally fits in the framework of topical
problems of Connes's noncommutative geometry~\cite{Con1}.

\subsubsection*{Aim of the paper}

This paper deals with elliptic theory on manifolds with corners.

Operators on manifolds with corners have been actively studied, and a number of
important interesting results emerged recently. For example, the $C^*$-closure
of symbol algebras was studied in~\cite{MeNi2}, and a spectral sequence
converging to the $K$-theory of the $C^*$-algebra of symbols was constructed.
Monthubert~\cite{Mon3} obtained a description of the operator algebra in the
spirit of noncommutative geometry in terms of a special groupoid that can be
associated with a manifold with corners (see also \cite{LeMo1}).
Bunke~\cite{Bun2} constructed index theory of Atiyah--Patodi--Singer type for
Dirac operators and studied cohomological obstructions to elliptic problems
(see also \cite{Loy1,LaMo3}); Monthubert and Nistor~\cite{MoNi1} produced a
formula for the boundary map in the $K$-theory of symbol algebras in
topological terms. Krainer~\cite{Kra1} studied boundary value problems in this
setting.

Although these results permitted finding the group classifying the homotopy
classes of elliptic operators in a number of special cases (e.g.,
see~\cite{MePi2} or \cite{Nis1}), the homotopy classification problem remained
open.

We solve Gelfand's problem for \emph{manifolds with corners}. Our goal is
to obtain a simple explicit formula for the classifying group in terms of
Atiyah's $K$-homology functor \cite{Ati4}.

\subsubsection*{Elliptic operators and $K$-homology}

Note that the idea of classifying elliptic operators by the $K$-homology
functor has long been known. For the reader's convenience, we recall it
using operators on a smooth compact manifold $M$ as an example.

The commutator of an elliptic zero-order\footnote{Working solely with
zero-order operators does not result in loss of generality, since order
reduction (say, multiplication by an appropriate power of the Laplace
operator) is always available.} operator $D$ on $M$ with the operator of
multiplication by a continuous function $f\in C(M)$ is compact
\begin{equation}\label{compa1}
[D,f]\in \mathcal{K}.
\end{equation}
By one definition, the \emph{contravariant $K$-theory} $K^0(C(M))$ of the
algebra $C(M)$ just consists of Fredholm operators for which the
commutators~\eqref{compa1} are compact. Thus $D$ determines an element of the
group $K^0(C(M))$, which is isomorphic to the $K$-homology group of $M$:
$$
K^0(C(M))\simeq K_0(M)
$$
by the Atiyah--Brown--Douglas--Fillmore--Kasparov theorem. Thus, assigning
the corresponding class in the $K$-homology to each elliptic operator, we
obtain a mapping
$$
\Ell(M)\longrightarrow K_0(M),
$$
where $\Ell(M)$ is the group of elliptic operators in sections of bundles
on $M$ modulo stable homotopy and  $K_0(M)$ is the even $K$-homology group
of ${M}$.

Kasparov showed this mapping to be an isomorphism. In other words, the
$K$-homology group of a smooth manifold classifies elliptic operators on
this manifold modulo stable homotopy.

This approach to classification also proved fruitful in the case of compact
stratified manifolds with singularities. Namely, it was shown
in~\cite{R:NaSaSt3} that in this case the even $K$-homology group of the
underlying compact topological space classifies elliptic operators on this
manifold.

However, no classification results were known for manifolds with corners of
codimension $\ge 2$. The classification in the form of the $K$-homology of the
manifold with corners, which suggests itself, is too meagre to be true: one can
always smooth the corners, and we see that the $K$-homology of the manifold
with corners is too coarse an invariant, for it does not take into account the
structure of a manifold with corners.

Moreover even the space whose $K$-homology would classify elliptic
operators was unknown.

\subsubsection*{Main result}
We establish the isomorphism
\begin{equation}
\label{funda1} \Ell(M)\simeq K_0({M}^\#),
\end{equation}
where $M$ is a manifold with corners and ${M}^\#$ is the dual manifold (see
\cite{Part1}) which is a stratified manifold with singularities. More
precisely, the isomorphism \eqref{funda1} will be established under the
following assumption concerning the combinatorial structure of the faces of
our manifold:
\begin{center}
\textit{The normal bundles of all faces of $M$ are trivial.}
\end{center}
If this assumption fails, then, generally speaking, the isomorphism
\eqref{funda1} does not hold. In this case, one should abandon the search for a
classifying \emph{space} and seek some \textit{algebra} whose $K$-cohomology
classifies elliptic operators. This algebra proves to be noncommutative, and
one needs to use ideas on noncommutative geometry. These results will be
considered elsewhere.

Note an interesting special case: if a manifold with corners is a
polyhedron with a given triangulation of the boundary, then the dual
stratified space is also a polyhedron, namely, the one used in the
classical proof of Poincar\'e duality in cohomology! For example, if $M$ is
a cube, then ${M}^\#$ is an octahedron. Thus the construction of the dual
manifold in the first part~\cite{Part1} of the present paper generalizes
the \textit{Poincar\'e dual polyhedron} to the case of noncontractible
faces.

\subsubsection*{Manifolds with corners and manifolds with multicylindrical ends}

\begin{figure}
\begin{center}
\includegraphics[height=5cm]{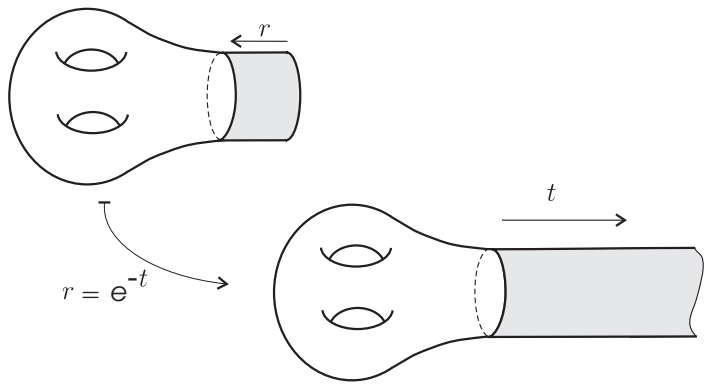}
\caption{Transition from a neighborhood of the
boundary to an infinite cylinder.}\label{ris4}
\end{center}
\end{figure}

Note that there is a different perspective on the theory of operators on
manifolds with corners. An application of a logarithmic change of variables
in a neighborhood of the boundary taking the boundary to infinity (see
Fig.~\ref{ris4}, where this is shown for manifolds with boundary) results
in the class of so-called \emph{manifolds with multicylindrical ends}.
These two pictures give the same operator algebras. Thus the results of the
present paper also provide classification on manifolds with
multicylindrical ends.

\subsubsection*{Outline of the paper}

This is the second of the two parts of the paper. In the first
part~\cite{Part1}, the dual manifold of a manifold with corners was
constructed and calculus of pseudodifferential operators ($\psi$DO) on
manifolds with corners was developed in the $C^*$-algebraic context.

The present part has the following structure.

In the first section, we recall some information from~\cite{Part1}. In
Sec.~2, we state the classification theorem. The proof occupies the next
three sections. Note that the general scheme of the proof is the same as
in~\cite{R:NaSaSt3}, and we proceed by induction, passing from a smooth
manifold to increasingly complex manifolds with singularities. In Sec.~6 we
discuss the relationship with some results due to Monthubert and Nistor. As
a consequence of the classification theorem, we obtain a $K$-homology
criterion for the vanishing of the index and a formula for the $K$-group of
the $C^*$-algebra of $\psi$DO with zero interior symbol (this algebra
corresponds to the $C^*$-algebra of the groupoid constructed by
Monthubert). In the appendix, we prove a higher analog of the relative
index theorem, which naturally arises when we obtain the classification of
operators.

\section{Manifolds with Corners and Dual Stratified Manifolds}\label{pargeom}

\subsubsection*{Manifolds with corners and faces}

Here we recall some information given in~\cite{Part1}.

Consider a manifold $M$ with corners of depth $k$. It has a natural
stratification
$$
M=\bigcup_{j=0}^{k} M_j,
$$
where the stratum $M_0=\ciM$ is just the interior of $M$ and each stratum
$M_j$ is the union of connected components, open faces $M_{j\a}$ of
codimension $j$ in $M$.

Each face $F=M_{ja}$ in the stratum $M_j$ is isomorphic to the interior of a
manifold $\overline{F}=\ov{M}_{ja}$, which will be called a closed face of $M$.
Faces of codimension one are called \emph{hyperfaces}.

\subsubsection*{Main assumption}

The main results of the paper will be obtained under the following
assumption.
\begin{assumption}\label{assa1}
The normal bundle $N_+F$ of each face $F$ is trivial.
\end{assumption}
In this case, the local defining functions $\rho_1,\ldots,\rho_j$ of $F$ are
globally defined as functions on the normal bundle $N_+F$.

\begin{remark*}
Assumption \ref{assa1} holds if all hyperfaces are \emph{embedded}, i.e., if
there exists a global defining function for each hyperface $F\subset M$.
However, it also holds, for some manifolds with nonembedded hyperfaces, say,
for the raindrop. The simplest example of a manifold with corners that does not
satisfy Assumption~\ref{assa1} is the Klein bottle with raindrop instead of the
circle as the base.
\end{remark*}

\subsubsection*{The dual space}
The dual space $M^\#$ of a manifold $M$ with corners was introduced
in~\cite{Part1}. If the original manifold is represented as the union
$$
M=\bigcup_{j\ge 0} M_j,\qquad M_j=\bigcup_{\alpha} M_{j\alpha},
$$
then ${M}^\#$ is the union of  \emph{dual faces},
$$
{M}^\#=\bigcup_{j\ge 0} {M}^\#_j,\qquad {M}^\#_j=\bigcup_{\alpha}
{M}^\#_{j\alpha},
$$
each of which is isomorphic to the interior of a simplex
$$
 {M}^\#_{j\alpha}\simeq \Delta^\circ_{j-1}.
$$
Here, by definition, ${M}^\#_0:=M_0$ is the interior of $M$. Thus to each
face $F$ of codimension $j$ in $M$ there corresponds a simplex $F^\#$ of
dimension $j-1$ in the dual space.

\subsubsection*{The fibration structure on $M^\#$}

It was proved in~\cite{Part1} that a neighborhood ${U}^\#$ of the stratum
${F}^\#$ is homeomorphic to the product of ${F}^\#$ by the cone
$$
K_{\Omega}=[0,1)\times \Omega/\{0\}\times \Omega
$$
whose base $\Omega$ is the dual space ${\overline{F}}^\#$ of the closed
face $\overline{F}$ (which is well defined, since $\overline{F}$ itself is
a manifold with corners). As a result, we find that ${M}^\#$ is a
stratified manifold with singularities.

\section{Classification Theorem}\label{maina}

Let $M$ be a manifold with corners satisfying Assumption \ref{assa1}, and
let $\Psi(M)$ be the $C^*$-algebra of zero-order $\psi$DO in the space
$L^2(M)$ (see~\cite{Part1}). The notion of a $\psi$DO acting on sections of
finite-dimensional vector bundles on $M$ is introduced in the usual way.
There is a natural equivalence relation, \textit{stable homotopy}, on the
set of elliptic operators. Recall the definition.
\begin{definition}
Two elliptic operators
$$
D:L^2(M,E)\to L^2(M,F)\;\;\text{ and }\;\; D':L^2(M,E')\to L^2(M,F')
$$
on $M$ are said to be \textit{stably homotopic} if there exists a
continuous homotopy
$$
D\oplus 1_{E_0}\sim f^*\bigl(D'\oplus 1_{F_0}\bigr)e^*
$$
of elliptic operators, where $E_0,F_0\in\Vect(M)$ are vector bundles and
$$
e:E\oplus E_0\longrightarrow E'\oplus F_0,\qquad f:F'\oplus
F_0\longrightarrow F\oplus E_0
$$
are bundle isomorphisms.
\end{definition}
Here ellipticity is understood as the invertibility of all components of
the symbol of the operator, and only homotopies of $\psi$DO preserving
ellipticity are considered.

\subsubsection*{Even groups $\Ell_0({M})$}
Stable homotopy is an equivalence relation on the set of elliptic $\psi$DO
acting in sections of vector bundles. By $\Ell_0({M})$ we denote the
corresponding quotient set. It is a group with respect to the direct sum of
operators, and the inverse in this group is given by the coset of the
almost inverse operator (i.e., an inverse modulo compact operators).

\subsubsection*{Odd groups $\Ell_1({M})$.}
Odd elliptic theory $\Ell_1({M})$ is defined in a similar way as the group
of stable homotopy classes of elliptic self-adjoint operators.
Stabilization is defined in terms of the operators $\pm Id$.

\begin{remark}
An equivalent definition of the odd $\Ell$-group can be given in terms of
smooth operator families on ${M}$ with parameter space $\mathbb{S}^1$
modulo constant families.
\end{remark}

We shall compute the groups $\Ell_*(M)$ for  $*=0$ and $*=1$, i.e., find
the classification of elliptic operators modulo stable homotopy. Our
approach is based on the following fact (see the definition of $\psi$DO
in~\cite{Part1}):
\begin{center}
\emph{$\psi$DO on $M$ can be viewed as local operators in the sense of
Atiyah on the dual manifold ${M}^\#$.}
\end{center}
Thus an elliptic $\psi$DO defines a Fredholm module on the space
$L^2(M)$ viewed as a $C({M}^\#)$-module. (For Fredholm modules and
$K$-theory, see~\cite{HiRo1} or \cite{Bla1}).

\subsubsection*{Classification of elliptic operators}

The following theorem is the main result of this paper.
\begin{theorem}\label{th1}
The mapping that takes each elliptic $\psi$DO to the corresponding Fredholm
module defines the group isomorphism
\begin{equation}\label{iso1}
\Ell_*(M)\stackrel\simeq\longrightarrow K_*({M}^\#).
\end{equation}
\end{theorem}
We shall obtain this theorem as a special case of the following more
general theorem.

\subsubsection*{Classification of partially elliptic operators \rom(cf.~\cite{R:NaSaSt3}\rom).}
Let $\Ell_*\left(M_{\ge j}\right)$ be the group generated by operators
whose symbols are invertible on the main stratum and all faces of
codimension $\ge j$. Thus we consider operators satisfying the ellipticity
condition on part of the faces.

The corresponding dual space
$$
{M}^\#_{\ge j}:={M}^\#\setminus \bigcup^{j-1}_{j'=1}{M}^\#_{j'}
$$
is obtained from ${M}^\#$ by
deleting all simplices of dimension $\le j-2$.
\begin{lemma}\label{lem1a}
An operator $D$ such that $[D]\in \Ell_*\left(M_{\ge j}\right)$ defines a
Fredholm module over the algebra $C_0({M}^\#_{\ge j})$ of functions on the
dual space.
\end{lemma}
\begin{proof}
We should verify the following properties of a Fredholm module: the
expression
$$
 f(DD^*-1)
$$
is compact for all $f\in C_0\left({M}^\#_{\ge j}\right)$ (here we assume
that $D$ is normalized by the condition $\sigma^*_s(D)=\sigma_s(D)^{-1}$
for $s\ge j$). The compactness follows from the fact that, by construction,
on each face $F\subset M$ either the corresponding symbol of our operator
is invertible  or the functions in the algebra
$C_0\left({M}^\#_{\ge j}\right)$ are zero.
\end{proof}
By Lemma~\ref{lem1a}, the mapping
\begin{equation}\label{iso2}
\Ell_*\left(M_{\ge j}\right)\stackrel{\varphi_j}\longrightarrow
K_*\left({M}^\#_{\ge j}\right)
\end{equation}
that takes partially elliptic operators to the corresponding Fredholm
modules is well defined (cf.~\cite{R:NaSaSt3}).

\begin{theorem}\label{th2}
For each $1\le j\le k+1$, the mapping \eqref{iso2} is an isomorphism.
\end{theorem}
Theorem \ref{th1} is the special case of Theorem \ref{th2} for $j=1$.

\section{Beginning of Proof of the Classification Theorem}
We prove Theorem~\ref{th2} by induction on $j$ decreasing from $k+1$ (where
$k$ is the depth of $M$) to~$1$.

\subsubsection*{Inductive assumption}
For $j=k+1$, the group $\Ell_*(M_{\ge j})$ classifies elliptic interior
symbols and hence is isomorphic to $K_c^*(T^*M)$. Moreover, the mapping
taking the symbol to the corresponding operator determines an isomorphism
$K_c^*(T^*M)\simeq K_*(M_0)$ (e.g., see~\cite{Kas3}). On the other hand,
the right-hand side of \eqref{iso2} in this case just contains the group
$K_*(M_0)$. Thus the theorem holds for $j=k+1$.

\subsubsection*{Inductive step.}
To justify the inductive step, we need to study exact sequences in
$K$-homology and $K$-theory permitting one to relate the maps $\varphi_j$
in \eqref{iso2} for two values of the subscript, $j$ and $j+1$.

\subsection{Exact sequence in $K$-homology (see~\cite{R:NaSaSt3})}\label{geom}

Consider the embedding
\begin{equation}\label{pair1}
{M}^\#_{\ge j}\supset{M}^\#_{j} .
\end{equation}
The complement ${M}^\#_{\ge j}\setminus{M}^\#_{j}$ is obviously equal to
${M}^\#_{\ge j+1}$, and we have the exact sequence of the
pair~\eqref{pair1} in $K$-homology
\begin{equation}\label{seq1}
\ldots\rightarrow K_*({M}^\#_{j})\rightarrow K_*({M}^\#_{\ge j})\rightarrow
K_*({M}^\#_{\ge j+1})\stackrel\partial\rightarrow
K_{*+1}({M}^\#_{j})\rightarrow\ldots.
\end{equation}
All maps but the boundary map $\partial$ in this sequence correspond to a
change of module structure on the corresponding Fredholm modules. The
boundary map $\partial$ can be reduced to a  form convenient for
computations by the following standard method.

Let ${U}^\#\subset {M}^\#_{\ge j}$ be the open neighborhood of the stratum
${M}^\#_{j}$ constructed\footnote{The neighborhood $U$ is defined as the
union of neighborhoods of all simplices $F^\#\subset M^\#_j$. By
construction, these neighborhoods are disjoint.} in~\cite[Sec.~1.2]{Part1}.
We have the homeomorphism
\begin{equation}
{U}^\#\simeq {M}^\#_{j}\times K_{{\overline{M}_j}^\#},
\end{equation}
where the cone $K_{\Omega_j}$ is the disjoint union of the cones
corresponding to the connected components of the base $\Omega_j$.

Then we have the mappings
$$
{M}^\#_{j}\times (0,1)\stackrel{\pi}{\longleftarrow}{M}^\#_{j}\times
(0,1)\times \Omega_j\stackrel{\simeq}\leftarrow {U}^\#\setminus
{M}^\#_{j}\stackrel{l}\longrightarrow {M}^\#_{\ge j+1}
$$
(by $l$ we denote the embedding of an open manifold, and $\pi$ is the
projection onto the first two factors), which permit us to represent the
boundary map $\partial$ in \eqref{seq1} as the composition
\begin{equation}\label{Kgranica}
K_*({M}^\#_{\ge j+1})\stackrel{l^*}\rightarrow K_*({U}^\#\setminus
{M}^\#_{j})\stackrel{\pi_*}\longrightarrow K_*((0,1)\times
{M}^\#_{j})\stackrel\beta\simeq K_{*+1}({M}^\#_{j})
\end{equation}
of the restriction $l^*$ of operators to an open set, the push-forward
$\pi_*$, and the periodicity isomorphism $\beta$. This representation
follows from the fact that $\partial$ is natural.

\subsection{Exact sequence related to elliptic operators (see~\cite{MeNi2})}

Let $M$ be a manifold with corners of depth $k>0$, and let $j$, $1\le j\le
k$, be some number. We denote the algebra formed by the symbols
$(\sigma_0,\sigma_{j},\sigma_{j+1},\ldots,\sigma_k)$ of all $\psi$DO on $M$
by
$$
\Sigma_j=\im (\sigma_0,\sigma_{j},\sigma_{j+1},\ldots,\sigma_k).
$$
Then we have the short exact sequence of $C^*$-algebras
\begin{equation}\label{seq2}
0\to J\to \Sigma_j\to \Sigma_{j+1}\to 0.
\end{equation}
Here the ideal $J$ consists of the symbols
$(\sigma_0,\sigma_{j},\sigma_{j+1},\ldots,\sigma_k)$ in which all
components but $\sigma_{j}$ are zero. From the compactness criterion for
$\psi$DO and compatibility conditions for symbols (see~\cite{Part1}), we
see that under these conditions the symbol $\sigma_{j}$ is a tuple of
compact-valued families decaying at infinity, so that one has the
isomorphism
$$
J\simeq \bigoplus_{F\subset M_j} C_0(\mathbb{R}^{j},\mathcal{K}L^2(F)),
$$
where the sum is taken over faces $F$ of codimension $j$ in $M$.

By virtue of this isomorphism, we can write out the exact sequence in
$K$-theory corresponding to the short sequence \eqref{seq2} in the form
\begin{equation}\label{seq3}
\ldots\to K_*(J)\to K_*(\Sigma_j)\to
K_*(\Sigma_{j+1})\stackrel{\delta}\longrightarrow K_{*+1}(J)\to\ldots.
\end{equation}
Clearly,
$$
 K_*(J)\simeq K_*(C_0(\mathbb{R}^j))\oplus K_*(C_0(\mathbb{R}^j))\oplus\ldots=\mathbb{Z}^l,
$$
where $l$ is the number of connected components in $M_j$. In terms of this
isomorphism, the boundary map $\delta$ can be represented (for $*=1$) in
the following simple form. An arbitrary class
$$
 [\sigma]\in K_1(\Sigma_{j+1})
$$
is realized by an invertible symbol
$$
\sigma=(\sigma_0,\sigma_{j+1},\sigma_{j+2},\ldots,\sigma_k).
$$
(From now on, for brevity we carry out the computations only for $K$-theory
elements representable by scalar operators; the consideration of the matrix
case differs only in the awkwardness of formulas.) Take an arbitrary symbol
$\sigma_j$ compatible with $\si$. The symbol $\sigma_{j}$ defines an elliptic
family with parameters in $\mathbb{R}^{j}$, and the index of that family is
a well-defined element of the $K$-group with compact supports of the
parameter space. One has
$$
\delta[\sigma]=\ind \sigma_j\in \bigoplus_{F\subset
M_j}K_0(C_0(\mathbb{R}^j)).
$$
There is a similar expression for the boundary map for the case $*=0$. (To
obtain it, one can pass to the suspension.)

\subsection{Comparison of exact sequences}

Let us show that the sequences \eqref{seq1} and \eqref{seq3} can be
combined into the commutative diagram
\begin{equation}\label{maindiagram}
\begin{array}{rccccl}
\ldots\to K_*(J)\quad& \to K_*(\Sigma_j)\to & K_*(\Sigma_{j+1})&
\stackrel{\delta}\longrightarrow & K_{*+1}(J)& \to\ldots\\
\downarrow\varphi_0\quad & \downarrow\varphi_j & \downarrow\varphi_{j+1} & & \downarrow\varphi_0 \\
\ldots\rightarrow K_{*+1}({M}^\#_{j})& \to K_{*+1}({M}^\#_{\ge j
})\to & K_{*+1}({M}^\#_{\ge j+1})& \stackrel{\partial}\longrightarrow &
K_{*}({M}^\#_{j})& \to\ldots
\end{array}
\end{equation}
(The construction of this diagram and the verification of its commutativity
will be finished in Sec.~\ref{compare1}.)

\subsubsection*{First, we define the vertical maps in the diagram.}
Without loss of generality, we can assume that $M$ has no connected
components with empty boundary. (Everybody knows the classification on such
components.) Then for all $j$ we have the isomorphism
\cite{Sav8}\footnote{This isomorphism generalizes the well-known expansion
$K^1(S^*M)\simeq \Ell(M)\oplus K^1(M)$ on a smooth closed manifold $M$ on
which there exists a nonzero vector field. Elimination of closed components
permits us to claim that there exists a nonzero vector field in our
situation.}
\begin{equation}\label{dek1}
K_*(\Sigma_j)\simeq \Ell_{*+1}(M_{\ge j})\oplus K_*(C(M)).
\end{equation}
Hence we define the maps $\varphi_j$, $j\ge 1$, in diagram
\eqref{maindiagram} as the composition
$$
K_*(\Sigma_j)\longrightarrow \Ell_{*+1}(M_{\ge j})\longrightarrow
K_{*+1}({M}^\#_{\ge j })
$$
of the projection onto the $\Ell$-group and the quantization \eqref{iso2}.
Thus these maps are induced by quantization, which takes symbols to
operators.

It remains to define the map $\varphi_0$. Just as the other vertical arrows
in the diagram, it is defined by quantization, namely, by quantization of
symbols $\sigma=\sigma_{j}$ in the ideal $J$. The quantization of elements
of the ideal differs from the quantization of general elements of the
algebra $\Sigma_j$ only in that the operator is considered in the $L^2$
space in a small neighborhood $U$ of the stratum $M_{j}$ in $M$ constructed
in~\cite[Lemma~1.9]{Part1}. We denote the operator by $\widehat{\sigma}_j$.

Let us define a module structure on $L^2(U)$. To this end recall that
$U$ can be considered also as a subset of the positive quadrant $N'_+M_j$ of the logarithmic
normal bundle. Thus, this space is naturally a
$C_0({M}^\#_{j})$-module. (Elements of $C_0({M}^\#_{j})$ act on $N'_+M_j$ as
operators of multiplication by radially constant functions $f(y)$,
in logarithmic coordinates $y=-\ln
\rho$.) The verification of locality of operator with respect to this
module structure (i.e., proving that the operator $\widehat{\sigma}_j$
commutes with operators of multiplication by functions modulo compact operators)
is immediate, and hence for the element $[\sigma_j]\in K_{*+1}(J)$ we define the element
\begin{equation}\label{var1}
\varphi_0[\sigma_j]:=[\widehat{\sigma}_j]\in K_*({M}^\#_j).
\end{equation}

\subsubsection*{Diagram \eqref{maindiagram} commutes}

The commutativity of the middle square of the diagram follows directly from
definitions.

\begin{lemma}
The left square of diagram \eqref{maindiagram} commutes.
\end{lemma}
\begin{proof}
Indeed, consider the composition of mappings passing through the right
upper corner of the square
$$
\begin{array}{rccccl}
K_*(J)& \to & K_*(\Sigma_j)\\
\downarrow\varphi_0 & & \downarrow\varphi_j \\
K_{*+1}({M}^\#_{j})& \to & K_{*+1}({M}^\#_{\ge j }).
\end{array}
$$
It takes an elliptic symbol $\sigma_j$ to the operator
$\widehat{\sigma}_j$, which acts outside a neighborhood $U$ of the stratum
${M}^\#_j$ as the identity operator (modulo compact operators) in the space
$L^2(M)$ with the natural structure of a $C_0({M}^\#_{\ge j
})$-module. Now if we restrict the operator to a neighborhood of ${M}_j$
(the element in $K$-homology remains unchanged, since the corresponding
Fredholm module changes by a degenerate module) and then use a homotopy to
reduce the module structure to the composition
$$
C_0({M}^\#_{\ge j })\stackrel{i^*}\to C_0({M}^\#_j)\stackrel{\pi^*}\to
C({U}^\#)\to \mathcal{B}(L^2(U)),
$$
where $\pi:{U}^\#\to {M}^\#_j$ is a projection and $i:{M}^\#_j\subset
{M}^\#_{\ge j}$ is an embedding, then we obtain the Fredholm module
assigned to the symbol $\sigma_j$ by the composition of maps passing
through left bottom corner of the square. The commutativity of the square
is thereby established.
\end{proof}

Verification of the commutativity of the square containing the boundary
maps is rather cumbersome, and so we make it in a separate section.

\section{Boundary and Coboundary Maps}

In this section, we establish the commutativity of the square
\begin{equation}\label{kvad}
\begin{array}{ccc}
 K_*(\Sigma_{j+1})&
\stackrel{\delta}\longrightarrow & K_{*+1}(J)\\
 \downarrow\varphi_{j+1} & & \downarrow\varphi_0 \\
 K_{*+1}({M}^\#_{\ge j+1})& \stackrel{\partial}\longrightarrow &
K_{*}({M}^\#_{j})
\end{array}
\end{equation}
containing the boundary maps in diagram \eqref{maindiagram}. The scheme of
proof is as follows. We
\begin{enumerate}
    \item Compute the composition $\varphi_0\circ\delta$.
    \item Compute the composition $\partial\circ \varphi_{j+1}$.
    \item Compare the resulting expressions.
\end{enumerate}

\subsection{Composition $\varphi_0\circ\delta$}

Let $[\sigma]\in K_*(\Sigma_{j+1})$ be the element defined by some symbol
$\sigma=(\sigma_0,\sigma_{j+1},\ldots,\sigma_k)$. Take a symbol $\sigma_j$
on $M_j$ compatible with $\si$ and denote by
\begin{equation}\label{sij}
 \widehat{\sigma}_j:L^2(NM_j)\to L^2(NM_j)
\end{equation}
the corresponding translation-invariant infinitesimal operator. (It is
conjugate to $\sigma_j$ by the Fourier transform.)

Representing the space $N'_+M_j$ as the product $ N'_+M_j\simeq M_j\times
{\mathbb{R}_+^j}$, we see that $L^2(NM_j)$ is a
$C_0(\bigsqcup\mathbb{R}^j_+)$-module. Here $\bigsqcup\mathbb{R}^j_+$ is
the disjoint union of as many open quadrants as there are faces of
codimension $j$ in $M$. The operator $\widehat{\sigma}_j$ is local with
respect to this module structure. We denote the corresponding element of
the $K$-homology group by
\begin{equation}\label{var2}
[\widehat{\sigma}_j]\in K_{*+1}(\bigsqcup\mathbb{R}^j_+).
\end{equation}

\begin{lemma}\label{l1}
The element $[\sigma]\in K_*(\Sigma_{j+1})$ satisfies the chain of
relations
\begin{equation}\label{lem2}
\varphi_0\delta[\sigma]=\varphi_0(\ind
\sigma_j)=\beta[\widehat{\sigma}_j]\in K_*({M}^\#_{j}),
\end{equation}
where
$\beta:K_{*+1}(\bigsqcup\mathbb{R}_+^j)=K_{*+1}({M}^\#_{j}\times\mathbb{R}_+)
\longrightarrow K_*({M}^\#_{j})$ is the Bott periodicity isomorphism, and
the index is understood as the index
$$
\ind \sigma_j\in K^{*+1}(\bigsqcup\mathbb{R}^j)\simeq K_{*+1}(J)
$$
of the elliptic operator-valued symbol $\sigma_j$.
\end{lemma}
\begin{proof}
For brevity, we assume that $M_j$ consists of exactly one stratum. In this
case, we have $M_j^\#\times\mathbb{R}_+\simeq \mathbb{R}_+^j$.

The first relation in \eqref{lem2} follows from definitions (since the
boundary map in $K$-theory of algebras is the index map).

1. Let us establish the second relation $\varphi_0(\ind
\sigma_j)=\beta[\widehat{\sigma}_j]$. The proof is based on the diagram
\begin{equation}\label{triang1}
 \xymatrix{
  K_{*+1}(\Sigma_{j+1})\ar[dr]
  \ar[r]^{\ind\sigma_j}& K_*(J)\ar[rd]^{\varphi_0}\ar[d]_q\\
  &
  K_*({M}^\#_j\times \mathbb{R}_+)\ar[r]^\beta &
  K_{*+1}({M}^\#_j),
 }
\end{equation}
where the map $K_{*+1}(\Sigma_{j+1})\to K_*({M}^\#_j\times \mathbb{R}_+)$
is induced by the map that takes the symbol $\sigma$ to the operator
$\widehat{\sigma}_j$ in \eqref{var2}. Finally, the group $K_*(J)\simeq
K^*(\mathbb{R}^j)$ is interpreted as the $K$-group of the cotangent bundle
to $\mathbb{R}^j_+$, and the map $q$ is induced by standard quantization
(to a symbol on the cosphere bundle, one assigns a pseudodifferential operator).

2. We claim that diagram~\eqref{triang1} commutes. Indeed, let us verify
the commutativity of the left triangle, i.e., the relation
\begin{equation}\label{luka3}
[\widehat{\sigma}_j]=q[\ind \sigma_j].
\end{equation}
Note  that the operator $\widehat{\sigma}_j$ is given over the product
$NM_j=M_j\times \mathbb{R}^j$. Moreover, it can be viewed as a $\psi$DO on
$\mathbb{R}^j$ with operator-valued symbol $\sigma_j=\sigma_j(\xi),$
$\xi\in \mathbb{R}^j$. This symbol is independent of the physical variables
$x\in \mathbb{R}^j$. Without loss of generality, it can be assumed to be
smooth with respect to the parameter $\xi$ (since $\sigma_j(\xi)$, just as
any $\psi$DO with a parameter, can be arbitrarily closely approximated by a
smooth $\psi$DO with a parameter; see~\cite{Part1}). Hence $\sigma_j(\xi)$
is an operator-valued symbol in the sense of \cite{Luk1}, i.e., has a
compact variation with respect to $\xi$ and all of its derivatives starting
from the first decay at infinity. Now relation \eqref{luka3} follows by
analogy with the generalized Luke theorem in~\cite{R:NaSaSt3}.

The commutativity of the right triangle follows  (see Corollary~\ref{aux1}
in the appendix) from the higher relative index theorem.

3. The commutativity of diagram~\eqref{triang1} implies the second relation
\eqref{lem2}. (The right-hand side is obtained if from the left top corner
of the diagram we go directly to the group $K_*(M_j^\#\times\mathbb{R}_+)$
and then apply the periodicity isomorphism $\beta$.)
\end{proof}

\subsection{Composition $\partial\circ\varphi_{j+1}$}

\subsubsection*{Space $N'_+M_j$ as a manifold with corners}

The image of the positive quadrant $N'_+M_j$ under the inverse of the
logarithmic map is the set $M_j\times [0,1)^j\subset M_j\times
\overline{\mathbb{R}}^j_+=N_+M_j$. Hence we treat $N'_+M_j$ as the interior
of a manifold with corners, denoted by $\overline{N'_+M_j}$. We denote the
corresponding dual space by ${\overline{N'_+M_j}}^\#$. On the complement
${\overline{N'_+M_j}}^\#\setminus {M}^\#_j$, there is a well-defined
projection
\begin{equation}\label{proek}
\begin{array}{ccc}
{\overline{N'_+M_j}}^\# \setminus {M}^\#_j& \stackrel{\pi}\longrightarrow &
{M}^\#_j \times\mathbb{R}_+\\
\displaystyle(y,x,\omega) &\mapsto &\left(\displaystyle\frac{y}{|y|},
\frac{|x|+1}{|y|}\right),
\end{array}
\end{equation}
whose fiber is the space ${\overline{M}_j}^\#$.

\subsubsection*{Reduction into a neighborhood of the edge}
We have the diagram of embeddings
\begin{equation}
\label{emb}
\begin{array}{rcl}
{M}^\#_j\subset& {U}^\# & \subset {M}^\#_{\ge j}\\
& \bigcap &\\
& {\overline{N'_+M}_j}^\#. &
\end{array}
\end{equation}

Let $[\sigma]\in K_{*+1}(\Sigma_{j+1})$ be the element determined by the
symbol $\sigma$ as above. By passing to the corresponding operators, we
obtain the element
$$
\varphi_{j+1}[\sigma]\in K_*({M}^\#_{\ge j+1}).
$$
On the other hand, the infinitesimal operator $\widehat{\sigma}_j$
compatible with $\sigma$ (see~\eqref{sij}) defines the element
$$
[\widehat{\sigma}_j]'\in K_*({\overline{N_+M}_j}^\#\setminus {M}^\#_j).
$$
This element is well defined, since the components of its symbol are
elliptic on the corresponding strata. We use primes to distinguish this
element from the element \eqref{var2}: although they are determined by one
and the same operator, the module structures on the $L^2$-spaces are different.

The naturality of the boundary map in $K$-homology results in the following
lemma.
\begin{lemma}\label{l2}
One has
$$
\partial \varphi_{j+1}[\sigma]=\partial'[\widehat{\sigma}_j]',
$$
where $\partial':K_*({\overline{N_+M}_j}^\#\setminus {M}^\#_j)\to
K_{*+1}({M}^\#_j)$ is the boundary map for the pair
${M}^\#_j\subset{\overline{N_+M}_j}^\#$.
\end{lemma}
\begin{proof}
Let $D$ be some operator on $M$ with symbol $\sigma$.

1. The infinitesimal operator $\widehat{\sigma}_j$ is obtained from $D$ by
localization to the set $M_j^\#$. Hence the restrictions $D|_U$ and
$\widehat{\sigma}_j|_U$ of these operators to a small space $U$ of $M_j$
are connected by a linear homotopy; i.e., one has
\begin{equation}\label{easy1}
[D|_{U}]=[\widehat{\sigma}_j|_U]\in K_*({U}^\#).
\end{equation}

2. By applying the naturality of the boundary map in $K$-homology to the
embedding diagram \eqref{emb}, we obtain
$$
\partial \varphi_{j+1}[\sigma]\equiv\partial [D]=\partial''[D|_{U}],
$$
where $\partial''$ is the boundary map for the pair ${M}^\#_j\subset
{U}^\#$. Now if on the right-hand side of the last relation we replace the
element $[D|_{U}]$ according to \eqref{easy1} and once more use the
naturality of the boundary map, then we obtain the desired relation
$$
\partial \varphi_{j+1}[\sigma]=\partial''[\widehat{\sigma}_j|_U]=
\partial'[\widehat{\sigma}_j]'.
$$
\end{proof}

Thus in what follows, when computing the composition $\partial\circ
\varphi_{j+1}$, we can (and will) work with the operator
$\widehat{\sigma}_j$ on $N'_+M_j$.

\subsubsection*{Homotopy of the module structure}

By \eqref{Kgranica}, the boundary map $\partial'$ in Lemma~\ref{l2} can be
represented as the composition
\begin{equation}\label{kompa1}
K_*({\overline{N'_+M}_j}^\#\setminus {M}^\#_j)
\stackrel{\pi_*}\longrightarrow
K_*({M}^\#_j\times\mathbb{R}_+)\stackrel\beta\rightarrow K_{*+1}({M}^\#_j)
\end{equation}
of the push-forward with respect to the projection $\pi$ and the
periodicity isomorphism.

Unfortunately, although the classes $[\widehat{\sigma}_j]$ and
$\pi_*[\widehat{\sigma}_j]'$ are determined by the same operator
$\widehat{\sigma}_j$, they have different module structures on the space
$L^2(NM_j)$: in the first case, the structure is independent of the
coordinate $x$, while in the second case it depends on
(see~\eqref{kompa1}).

Let us make a homotopy of module structures. To this end, we define a
homotopy
$$
 \pi^\varepsilon:{\overline{N'_+M}_j}^\#\setminus
{M}^\#_j\longrightarrow {M}^\#_j\times\mathbb{R}_+
$$
of projections by the formula (cf.~\eqref{proek})
$$
\pi^\varepsilon(y,x,\omega):=\left(\displaystyle\frac{y}{|y|},
\frac{\varepsilon|x|+1}{|y|}\right).
$$
This formula defines a continuous family of maps for $\varepsilon>0$.
However, the family is not continuous as $\varepsilon\to 0$.\footnote{And
hence the map $\pi^0_*$ is not defined on the $K$-group.} Nevertheless,
continuity takes place for the Fredholm modules, as shown by the following
lemma.
\begin{lemma}\label{l3}
The family $\pi^\varepsilon_*(\widehat{\sigma}_j)'$ of Fredholm modules
obtained by the change of module structure defines a homotopy in the sense
of $KK$-theory, and one has
\begin{equation}
\lim_{\varepsilon\to 0}\pi_*^\varepsilon (\widehat{\sigma}_j)'=\pi_*^0
(\widehat{\sigma}_j)',
\end{equation}
whence it follows that $\pi_*[\widehat{\sigma}_j]'=[\pi_*^0
(\widehat{\sigma}_j)]'\in K_*(M_j^\#\times\mathbb{R}_+)$.
\end{lemma}
\begin{proof}
For brevity, we assume that $M_j$ consists of a single face, i.e., is
connected. Then the homotopy in the sense of $KK$-theory means (e.g., see
\cite{Bla1}) that for each function $f\in C_0(\mathbb{R}^j_+)$ the family
$$
g^\varepsilon=(\pi^\varepsilon)^*(f):L^2(N'_+M_j)\longrightarrow
L^2(N'_+M_j)
$$
of operators of multiplication by the functions $(\pi^\varepsilon)^*(f)$ is
strongly $*$-continuous and that the operator families
$$
[g^\varepsilon,\widehat{\sigma}_j],\quad
g^\varepsilon(\widehat{\sigma}_j\widehat{\sigma}_j^{-1}-1)
$$
in $L^2(NM_j)$ are continuous families of compact operators as
$\varepsilon\to 0$.

It suffices to prove all these facts for (a dense set of) smooth functions
$f$. If $f$ is smooth, then one should smooth the family $g^\varepsilon$ and
use the composition formulas, which provide the desired compactness and
continuity.
\end{proof}

\subsection{Comparison of the compositions $\varphi_0\circ\delta$ and
$\partial\circ\varphi_{j+1}$}\label{compare1}
 Now let us use Lemmas \ref{l1}--\ref{l3}.
We obtain the chain of relations
\begin{multline*}
\partial\varphi_{j+1}[\sigma]\stackrel{\text{Lemma \ref{l2}}}=
\partial'[\widehat{\sigma}_j]'\stackrel{\text{formula \eqref{kompa1}}}
=\beta \pi^1_*[\widehat{\sigma}_j]' \stackrel{\text{Lemma \ref{l3}}}=\beta
[\pi^0_*\widehat{\sigma}_j]'=\\
=\beta[\widehat{\sigma}_j] \stackrel{\text{Lemma
\ref{l1}}}=\varphi_0\delta[\sigma].
\end{multline*}
The equality at the end of the first row corresponds to the identical
coincidence of the corresponding Fredholm modules.

Thus the square \eqref{kvad} commutes, and we have established the
commutativity of diagram \eqref{maindiagram}.

\section{End of Proof of the Classification Theorem}

By virtue of the isomorphism \eqref{dek1}, we can single out and cancel the
summand $K_*(C(M))$ in diagram \eqref{maindiagram} in the terms
$K_*(\Sigma_j)$ and $K_*(\Sigma_{j+1})$. We obtain the diagram
\begin{equation}\label{maindiagram1}
\begin{array}{rccccl}
\ldots\to K_*(J)& \to \Ell_{*+1}(M_{\ge j})\to & \Ell_{*+1}(M_{\ge j+1})&
\stackrel{\delta}\rightarrow & K_{*+1}(J)& \to\ldots\\
\downarrow\varphi_0\quad & \downarrow\varphi_j & \downarrow\varphi_{j+1} & & \downarrow\varphi_0 \\
\ldots\rightarrow K_{*+1}({M}^\#_{j})& \to K_{*+1}({M}^\#_{\ge j
})\to & K_{*+1}({M}^\#_{\ge j+1})& \stackrel{\partial}\rightarrow &
K_{*}({M}^\#_{j})& \to\ldots.
\end{array}
\end{equation}
The map $\varphi_{j+1}$ is an isomorphism by the inductive assumption. The
map $\varphi_0$ is also an isomorphism (see Corollary~\ref{aux1} in the
Appendix). Since the diagram commutes, we can apply the $5$-lemma and
obtain the desired justification of the induction step in Theorem
\ref{th2}: if the map $\varphi_{j+1}$ is isomorphic on the  $\Ell$-group,
then so is the map $\varphi_{j}$.

The proof of Theorem \ref{th2} is complete.

\section{Application to the Monthubert--Nistor index}

Let us discuss the relationship with the problems considered by Monthubert
and Nistor \cite{MoNi1}. In the notation of the present paper, for the case
of manifolds with embedded corners they considered the short exact sequence
\begin{equation}\label{shorts}
0\to J\longrightarrow \Psi(M) \stackrel{\sigma_0}\longrightarrow C(S^*M) \to 0,
\end{equation}
where $\sigma_0$ is the interior symbol map, and the ideal $J$ consists of
operators with zero interior symbol. They studied the boundary map
corresponding to this sequence:
$$
\delta:K_*(C(S^*M))\longrightarrow K_{*+1}(J).
$$
For a closed manifold $J$ is the ideal of compact operators (hence
$K_*(J)\simeq \mathbb{Z}$) and the boundary map coincides with the analytic
index. Moreover, Monthubert and Nistor showed that in the general case this map
has an important topological meaning: it gives the obstruction to the existence
of an invertible operator with a given interior symbol. For these reasons,
Monthubert and Nistor call this map the \emph{analytic index of manifolds with
corners}.

We claim that the classification theorem readily implies a $K$-homology
criterion for the vanishing of the analytic index. Indeed, consider the diagram
\begin{equation}
\label{mont1}
\begin{array}{cccl}
K_{*+1}({M}^\#) & \longrightarrow&
K_{*+1}(M_0)&\stackrel\partial\longrightarrow
K_{*}({M}^\#\setminus M_0)\\
\varphi_1\uparrow& &\uparrow\varphi_{k+1}\\
K_*(\Psi(M))&\longrightarrow& K_*(C(S^*M))&\stackrel\delta\longrightarrow
K_{*+1}(J),
\end{array}
\end{equation}
where the lower row is the sequence induced by the short exact sequence
\eqref{shorts} and the upper row is the exact sequence of the pair
${M}^\#\setminus M_0\subset {M}^\#$ in $K$-homology. The maps $\varphi_1$ and
$\varphi_{k+1}$ are induced by quantization of elliptic symbols on $M^\#$ and
$M_0$ correspondingly (cf.~\eqref{maindiagram}). The diagram is obviously
commutative.

From the exactness of the sequences and the obvious commutativity of the
diagram, we obtain the following assertion. Let us assume for simplicity that
$M$ has no connected components with empty boundary.
\begin{proposition}
The analytic index $\delta(x)\in K_{*+1}(J)$ of $x\in K_*(C(S^*M))$ vanishes if
and only if $\partial\varphi_{k+1}(x)=0$.
\end{proposition}
\begin{proof}
1. There are splittings (cf.~\eqref{dek1})
$$
K_*(\Psi(M))\simeq \widetilde{\Ell}_{*+1}(M)\oplus K_{*}(C(M)),\quad
K_*(C(S^*M))\simeq \Ell_{*+1}(M_0)\oplus K_{*}(C(M)),
$$
where $\widetilde{\Ell}$ is the \emph{reduced} $\Ell$-group generated by
operators of index zero. Moreover, the direct summands $K_{*}(C(M))$ can be
cancelled in~\eqref{mont1}. This does not affect the boundary map. Hence, we
obtain the commutative diagram
\begin{equation}
\label{mont2}
\begin{array}{cccl}
\widetilde{K}_{*+1}({M}^\#) & \longrightarrow&
K_{*+1}(M_0)&\stackrel\partial\longrightarrow
\widetilde{K}_{*}({M}^\#\setminus M_0)\\
\varphi_1\uparrow& &\uparrow\varphi_{k+1}\\
K_*(\Psi(M))/K_*(C(M))&\longrightarrow&
K_*(C(S^*M))/K_*(C(M))&\stackrel\delta\longrightarrow K_{*+1}(J),
\end{array}
\end{equation}
where $\widetilde{K}_*$ is the reduced $K$-homology group generated by
operators of index zero.

3. By  the classification theorem, the quantization maps $\varphi$ in
\eqref{mont2} induce isomorphisms. Hence, the commutativity of the diagram
readily shows that vanishing of $\delta$ is equivalent to the vanishing of the
boundary map $\partial$ in $K$-homology.
\end{proof}
The reader can readily rewrite this formula in a more explicit form as a
condition on the interior symbol $\sigma_0$. There is also a cohomological
form of this condition. Needless to say, the cohomological formula is only
valid modulo torsion.

\begin{remark*}
One actually has the group isomorphism
$$
K_*(J)\stackrel\simeq\longrightarrow \widetilde{K}_*({M}^\#\setminus M_0)
$$
determined by quantization of operators with zero interior symbol.  (One can
readily obtain this isomorphism by reproducing the proof of our classification
theorem. In the proof, only the inductive assumption is changed: now for
$j=k+1$ we claim that $0=0$.)
\end{remark*}

\appendix

\section{Higher Relative Index Theorem}

Consider the map (see Sec.~\ref{compare1})
$$
 \varphi_0:K_j(C_0(\mathbb{R}^j))\to
K_{j-1}(\Delta^\circ_{j-1}),
$$
induced by the map taking a symbol $\sigma(\xi),$  $\xi\in \mathbb{R}^j$,
to the corresponding translation-invariant operator
\begin{equation}\label{transl1}
\widehat{\sigma}:L^2(\mathbb{R}^j)\longrightarrow L^2(\mathbb{R}^j).
\end{equation}
Here the space $L^2(\mathbb{R}^j)$ is equipped with the following module
structure over the algebra $C_0(\Delta^\circ_{j-1})$ of functions on the
interior of the simplex: a function $f\in C_0(\Delta^\circ_{j-1})$ is
viewed as a radially constant function equal to zero outside the positive
quadrant  $\mathbb{R}_+^j$.

\begin{proposition}
For the index pairing of the element $\varphi_0[\sigma]\in
K_{j-1}(\Delta^\circ_{j-1})$ with an arbitrary element
$$
 [a]\in K_{j-1}(C_0(\Delta^\circ_{j-1}))\simeq \widetilde{K}^{j-1}(\mathbb{S}^{j-1}),
$$
where $\widetilde{K}$ is the reduced $K$-group, one has the formula
\begin{equation}\label{higher}
\langle\varphi_0[\sigma],a\rangle=\ind_t\Bigl([\sigma]\times [a]\Bigr),
\end{equation}
where the product $[\sigma]\times [a]$ is defined as the composition
$$
K^j(\mathbb{R}^j)\times \widetilde{K}^{j-1}(\mathbb{S}^{j-1})\to
K^1(S^*\mathbb{R}^j)\to K^0(T^*\mathbb{R}^j),
$$
and $\ind_t:K^0(T^*\mathbb{R}^j)\to \mathbb{Z}$ is the topological
Atiyah--Singer index for $\mathbb{R}^j$.
\end{proposition}
\begin{remark*}
For $j=1$, this assertion is reduced to the relative index theorem for
operators on manifolds with conical points. We mean the formula for the
difference of indices of operators with equal interior symbols, or,
equivalently, for the index of elliptic operators of the form $1+G$
$$
\ind (1+G)=w(1+g),
$$
where the interior symbol of $G$ is zero and $w(1+g)$ is the winding number
of the conormal symbol $1+g$. Hence the index formula~\eqref{higher} in the
general case can be referred to as the \emph{higher relative index
formula}.
\end{remark*}
\begin{proof}
To be definite, we consider the case of even $j$. (The odd case can be
considered in a similar way.)

1. The element
$$
[\sigma]\in K_0(C_0(\mathbb{R}^j))
$$
is determined by some projection-valued function $p(x)$ on $\mathbb{R}^j$
equal to the diagonal projection ${\rm diag}(1,0)$ at infinity. Conversely,
the element
$$
[a]\in K_{j-1}(C_0(\Delta^\circ_{j-1}))\simeq
\widetilde{K}^{j-1}(\mathbb{S}^{j-1})
$$
is determined by some invertible function $a(x)$ on the sphere
$\mathbb{S}^{j-1}$. To simplify the notation, we assume that this is a
scalar function. The matrix case can be considered in a similar way.

2. In this notation, the index pairing $\langle\varphi_0[\sigma],a\rangle$
is by definition equal to the index of the Toeplitz operator (e.g., see
\cite{HiRo1})
\begin{equation}\label{toep}
 \widehat{p}a:\im \widehat{p}\longrightarrow  \im \widehat{p},
\end{equation}
where $\widehat{p}:L^2(\mathbb{R}^j)\longrightarrow L^2(\mathbb{R}^j)$ ---
is the projection determined by the symbol $p$, as in \eqref{transl1}.

3. To compute the index of the operator \eqref{toep}, we make the Fourier
transform. Then the operator $\widehat{p}$ becomes the projection $p(x)$,
and the space $\im \widehat{p}$ becomes the space of sections of the bundle
given by the range of  $p(x)$. Conversely, the operator of multiplication
by $a$ passes into a translation-invariant $\psi$DO in $\mathbb{R}^j$ with
principal symbol $a=a(\xi)$. Hence we obtain
\begin{equation}\label{toep1}
 \langle\varphi_0[\sigma],a\rangle=\ind\left(
 p\widehat{a}:\im {p}\longrightarrow  \im {p}
 \right).
\end{equation}

4. The last operator coincides at infinity with a direct sum of the
invertible operator $\widehat{a}$ acting on functions. By the index
locality property, the difference of their indices is given by the
Atiyah--Singer formula, which gives the desired expression \eqref{higher}.
\end{proof}

Let us rewrite this index formula in the equivalent form used in the main
text.
\begin{corollary}\label{aux1}
The following triangle commutes\rom:
\begin{equation}\label{triang}
 \xymatrix{
  K_*(C_0(\mathbb{R}^j))\ar[rd]^{\varphi_0}\ar[d]_q\\
  K_*(\Delta^\circ_{j-1}\times (0,\infty))\ar[r]^\beta &
  K_{*+1}(\Delta^\circ_{j-1}),
 }
\end{equation}
where $\beta$ is the Bott periodicity isomorphism and
$$
q:K_*(C_0(\mathbb{R}^j))=K^*(T^*\mathbb{R}^j_+)\longrightarrow
K_*(\Delta^\circ_{j-1}\times (0,\infty))
$$
is the standard pseudodifferential quantization in
$\mathbb{R}^j_+=\Delta^\circ_{j-1}\times(0,\infty)$.
\end{corollary}
\begin{proof}
All groups in the triangle \eqref{triang} are isomorphic to $\mathbb{Z}$
and have natural Bott generators. It is known that the maps $q$ and $\beta$
take Bott elements to Bott elements. Hence to verify the commutativity of
the diagram it suffices to verify this property for $\varphi_0$. But this
readily follows from the index formula \eqref{higher}.
\end{proof}

\providecommand{\bysame}{\leavevmode\hbox to3em{\hrulefill}\thinspace}

\end{document}